\documentclass[11pt]{elsarticle}
\usepackage{amssymb}
\usepackage{amsmath}
\usepackage{graphicx}

\newcommand{\hamza}{\raisebox{.4ex}{$\rhook$}}

\title{Jost B\"{u}rgi's Method for Calculating Sines}

\author{Menso Folkerts, Dieter Launert, Andreas Thom}

\begin{document}

\maketitle

From various sources we know that the Swiss instrument maker and mathematician Jost B\"{u}rgi (1552-1632) found a new way of calculating any sine value. Many mathematicians and historians of mathematics have tried to reconstruct his so-called ``Kunstweg'' (``skillful / artful method''), but they did not succeed. Now a manuscript by B\"{u}rgi himself has been found which enables us to understand his procedure. The main purpose of this article is to explain B\"{u}rgi's method. It is totally different from the conventional way to calculate sine values which was used until the 17th century. First we will give a brief overview of the early history of trigonometry and the traditional methods for calculating sines.

\section{Historical remarks on trigonometry and on trigonometrical tables}

The main purpose of the following remarks is to show how values of trigonometric functions, especially of chords and sines, were calculated before the time of B\"{u}rgi. For this reason it is necessary to go back to Greek antiquity, the Arabic-Islamic tradition and the Western European Middle Ages. Of course this is not the place to present an extensive history of trigonometry%
\footnote{This can be found in [Van Brummelen, 2009].}.

In Greek antiquity, trigonometric functions were used in different contexts in order to calculate triangles and quadrangles: in geodesy particularly for determining heights and distances and in astronomy for calculating spherical triangles. In geodesy the main method was to use the proportion of the catheti in plane orthogonal triangles, i.e., in modern terms, the tangent. In astronomy the central problem was to find relations between the chords and the radius of the circle, today expressed by the sine.

When simple geodetic measurements had to be carried out, properties of similar triangles were used to find the fourth proportional with the help of the theorem of intersecting lines. For measuring angular distances on earth or in the sky the Jacob's staff was available.

Astronomical calculations were much more complicated. The main purpose was to solve problems of spherical geometry, i.e. to calculate spherical triangles and quadrangles. The Greeks had developed a method to calculate chords which was based on the fact that the central angle $\alpha$ of a circle is a measure for the circular arc (arc $\alpha$) and that there is a unique relation between the arc of a circle and the corresponding chord (crd $\alpha$). Accordingly they developed tables of chords as well as theorems on the relations between chords. In its fully developed form, this system is found in the \emph{Almagest}, the central astronomical work of Ptolemy (1st half of the 2nd c. AD). It was based upon the earlier achievements of Hipparchus (ca. 150 BC)%
\footnote{On Hipparchus, especially on his table of chords, see [Van Brummelen, 2009, 34-46].}.

The calculation of chords and its application to astronomical problems are taught in the first book of the Almagest. In Chapter 11 Ptolemy presents his table of chords. This is the oldest extant table of this kind. For any arc between 0$^{\circ}$ and 180$^{\circ}$ at intervals of half a degree it gives the value of the corresponding chord. The arcs and the chords are expressed in the sexagesimal system and the radius, too, is divided sexagesimally.

In chapter 10 Ptolemy explains in detail how he found these values%
\footnote{See [Van Brummelen, 2009, 70-77].}.
He proceeds as follows:

1. He starts with chords which correspond to the sides of special regular polygons. From Euclid's \emph{Elements} it was known that regular polygons with sides $n$ = 3, 4, 5, 6 and 10 can be constructed easily. The sides of these polygons are the chords of the angles 120$^{\circ}$, 90$^{\circ}$, 72$^{\circ}$, 60$^{\circ}$ and 36$^{\circ}$, accordingly. By theorems given in Euclid's \emph{Elements} Ptolemy calculates the chords of these five angles.

2. Then Ptolemy shows how to find the chord of the sum and difference of two arcs whose chords are known. His proof is based upon the theorem of the cyclic quadrilateral (the so-called ``Ptolemy's Theorem'').

3. Finally, Ptolemy gives a procedure on how to find the chord of half a given arc.

4. Now, starting with the known five chords, Ptolemy is able to compute the chord of $1\frac{1}{2}$ degrees%
\footnote{The steps are: crd 72$^{\circ}$, crd 60$^{\circ}$, crd 12$^{\circ}$, crd 6$^{\circ}$, crd 3$^{\circ}$, crd $1 \frac{1}{2}^{\circ}$.}
and from this the chords of all multiples of $1\frac{1}{2}^{\circ}$. But the chord of 1 degree cannot be found in this way. This has to do with the fact that the problem is closely related to the trisection of the angle. In order to find this value, Ptolemy uses a special interpolation procedure. It is based on the relation $\frac{\mathrm{crd}\; \alpha}{\mathrm{crd}\; \beta} < \frac{\alpha}{\beta}$, if $\beta < \alpha < 90^{\circ}$%
\footnote{See [Van Brummelen, 2009, 76f.].}.

After having found crd 1$^{\circ}$, Ptolemy can compute crd $\frac{1}{2}^{\circ}$ and, starting from this, the complete table of chords for any half degree.

Ptolemy's table of chords makes it possible to calculate angles and sides of arbitrary triangles and quadrangles. In the chapters of the \emph{Almagest} which follow the tables Ptolemy explains his methods, which are based upon the so-called ``Theorem of Menelaus''%
\footnote{See [Van Brummelen, 2009, 58f.].}.

Up to the 16th century, Ptolemy's methods and his table of chords dominated trigonometric calculations and theoretical astronomy. The \emph{Almagest} was translated from Greek into Arabic in the 9th century and from Arabic into Latin in the 12th century. Only in the 16th century its Greek version was recovered in Western Europe.

Trigonometric tables were developed not only in Greece, but independently in India, too. However, the Indians did not start from an arc, but used a value which we call today ``sine''. The sine of an angle $\varphi$ is half of the chord of the double angle $2\varphi$. This means:
$\sin \varphi = \frac{1}{2} \cdot \mathrm{crd}\; 2\varphi$.
Thus the sine of an angle and the chord of the double angle differ only by a constant factor, and therefore a chord table can be transferred easily into a sine table. But the substitution of the chord by the half-chord was important, because the sine is the side of an orthogonal triangle, so that it was possible to apply the theorems to orthogonal triangles, e.g. the Pythagorean Theorem, to the trigonometrical functions. The Indians did not use only the sine, but also the difference between the radius and the cosine (called ``versed sine''). Tables for these two functions were already available in India about 500 AD%
\footnote{On early sine tables in India and their calculation, see [Van Brummelen, 2009, 95-105].}. Indian trigonometry used a number of values for the radius of the base circle of the sine function (for example, $R=3438$), but never $R = 1$.

Through translations into Arabic, which were made in the 8th and 9th centuries, both the astronomical tables and writings of the Indians and Ptolemy's table of chords became known in the Arabic-Islamic world. Therefore in the Islamic countries both methods were known: the Greek concept of chords and the Indian concept of sine. In the following centuries, numerous collections of astronomical tables, the so-called \emph{Z\={\i}j}es, emerged. They served for the calculation of the movements of the heavenly bodies and for other astronomical and astrological purposes. Most \emph{Z\={\i}j}es contain sine tables and many of them also tables of the tangent.

It suffices to give some short remarks on sine tables in the Arabic-Islamic world%
\footnote{For a more extensive treatment, see [Van Brummelen, 2009, 135-222].}.
The oldest known collection of astronomical tables was arranged by Mu\d{h}ammad ibn M\=us\=a al-Khw\=arizm\=\i\ (around 820), but it has only survived in a Latin reworking. By this work Indian trigonometry was spread in the Arab countries. In the extant version of his \emph{Z\={\i}j al-Sindhind} the sine values are given to seconds for the radius $R = 60$, which is most likely of Greek origin. But other sources (e.g. Ibn al-Muthanna's commentary) indicate that Khw\=arizm\={\i}'s sine tables were originally based on the radius 150 which we often find in Indian texts.

From the end of the 9th to the 11th centuries some scholars computed more exact sine tables. Most important were al-Batt\=an\=\i, Ibn Y\=unus, Ab\=u'l-Waf\=a\hamza{} and al-B\={\i}r\=un\={\i}. In principle they used Ptolemy's method, but they applied more accurate interpolation methods in order to compute the value of sin 1$^{\circ}$%
\footnote{See [Van Brummelen, 2009, 137-147].}.

Of special importance in Europe were the ``Toledan tables'', a collection of astronomical tables which was compiled in the 11th century by a group of Andalusian scientists, the most prominent of whom was the astronomer al-Zarq\=allu. They used different sources, and in the course of time the work underwent some changes. We do not have the Arabic text, but only a large number of Latin manuscripts which represent different forms of transmission from the 12th century onwards. With the ``Toledan tables'' explanations (\emph{canones}) were provided. They, too, are extant in Latin only. Since Toomer's fundamental work%
\footnote{[Toomer, 1968]. The tables and their \emph{Canones} are now available in a complete critical edition [Pedersen, 2002].}
we know that the ``Toledan tables'' are closely connected with the tables of Khw\=arizm\={\i} and Batt\=an\={\i} and that most of the tables were taken over from either of these two works. Only two out of more than 80 tables are trigonometric: a sine table with radius 150 and another one with radius 60. Since the first table refers to the Indian value $R = 150$, it probably stems from Khw\=arizm\={\i}. The second table agrees in general with Batt\=an\={\i}'s \emph{Z\={\i}j}.

From the 12th century on, also in Europe Ptolemy's table of chords had to compete with other trigonometric tables which had been developed in India and had arrived via the Arabs to Western Europe. In the 13th and 14th centuries a wide variety of writings originated in the West which were based upon the \emph{Canones} of the ``Toledan tables''. Of special interest are the \emph{Canones tabularum primi mobilis} of John de Lineriis (1322)%
\footnote{Partial edition in [Curtze, 1900, 391-413].}
and the \emph{Quadripartitum} of Richard of Wallingford (before 1326)%
\footnote{Edited in [North, 1976, 1, 21-169].}.
Especially important were the \emph{Canones} of John de Lineriis. From the \emph{Canones} of the ``Toledan tables'' he adopted not only the method of calculating the sines, but also the sine table itself. It is calculated for every half degree to radius 60.

The methods which were used in the West for numeric calculation of sine values resemble Ptolemy's procedure. Most Western texts up to the 16th century include a figure which is already present in the \emph{Canones} to the ``Toledan tables''. It deals with finding the sine values of the so-called ``kardagae'', i.e. of multiples of 15$^{\circ}$%
\footnote{The word \emph{kardaga} refers to the Indian practice to divide the right angle into 24 parts. Later it was used for the division of the circle into 24 parts (see [Van Brummelen, 2009, 218-220]).}.
With the help of this figure all \emph{kardagae} can be found and the sine values can be computed easily. In order to determine further values, theorems are used to calculate the sine of half of an angle and the sine of the complementary angle. In this way one can compute the sine of 3$^{\circ}$ 45$'$. For the intermediate values interpolations are used which are, in general, similar to the Greek methods.

In the 15th century, scholars of the university in Vienna recalculated the sine values. The first of these was John of Gmunden (1380-1442). In his \emph{Tractatus de sinibus, chordis et arcubus} (1437) he gave two different methods: that from the ``Toledan tables'' and that from Ptolemy. As a matter of fact, he wrote two different treatises which led to the same result, i.e. to a sine table. John of Gmunden's work became known to Georg Peurbach (1421-1461) and to Johannes Regiomontanus (1436-1476). With a large amount of work Regiomontanus computed sine tables for every minute of arc. Three of these tables with increasing precision are extant: one for radius 60,000, the second for 6,000,000 and the third for 10,000,000. The third table is the first known table in which the sines are given decimally. It was finished in 1468, but published only in 1541. Different from the Indian and Islamic tables with fractional places, Regiomontanus calculated integer sine values and chose the radius of the base circle in order to obtain the precision he required%
\footnote{The contributions of these three mathematicians to trigonometry are summarized in [Van Brummelen, 2009, 248-263].}.

In the 16th century, the most important contributions to the calculation of trigonometric tables came from Nicolaus Copernicus (1473-1543), Georg Joachim Rheticus (1514-1574) and Bartholomaeus Pitiscus (1561-1613)%
\footnote{For details see [Van Brummelen, 2009, 265-282].}.
Copernicus had computed a sine table for his principal astronomical work by using Ptolemy's method. Rheticus arranged that the trigonometric part of \emph{De revolutionibus} was published in 1542, one year before Copernicus' main work appeared in print. Surprisingly, the edition which Rheticus initiated in 1542 does not contain the sine values calculated by Copernicus, but more exact values with 7 decimal places stemming from Regiomontanus. About ten years later, in 1551, Rheticus published his \emph{Canon doctrinae triangulorum}. This booklet brings something new: Rheticus related trigonometric functions directly to angles rather than to circular arcs, and he provided not only tables of the well-known four trigonometric functions (sin, cos, tan, cot), but also secant and cosecant tables. Obviously these values were computed by Rheticus himself. After 1551 he continued his calculations and around 1569 he had finished tables with an increment of 10 seconds. The work appeared in press only in 1596, after Rheticus' death, under the title \emph{Opus Palatinum}. Since the tables contained notable errors in some parts, a revision was soon considered to be necessary. This was accomplished by Pitiscus and was published in 1613 (\emph{Thesaurus mathematicus}). Pitiscus' publication marked a provisional end to the calculation of trigonometrical tables according to traditional methods.

\section{Jost B\"{u}rgi's method to compute the sine values}

All procedures for calculating chords and sines described so far are based in principle on the method which Ptolemy had presented in his \emph{Almagest}. Totally different from this is a procedure which Jost B\"{u}rgi (1552-1632) invented in the second half of the 16th century. With his method he was able to compute the sine of each angle with any desired accuracy in a relatively short time. B\"{u}rgi explained his procedure in a work entitled \emph{Fundamentum Astronomiae}. Up to now, this treatise, as well as his method of calculation, have been totally unknown. This is the first publication of B\"{u}rgi's procedure%
\footnote{A note about the treatise, but not about B\"{u}rgi's method of calculating sines, has been given in [Folkerts, 2014].}.

Jost B\"{u}rgi is well-known in the history of mathematics because he invented logarithms independently of John Napier%
\footnote{It seems that B\"{u}rgi invented the logarithms in the 1580s, but he published his table of logarithms not before 1620.}.
Another treatise of B\"{u}rgi is less known, since it was edited only in modern times, namely his algebraic work \emph{Coss}%
\footnote{Edited in [List, Bialas, 1973].}.
Unlike most other writings bearing this title, B\"{u}rgi's \emph{Coss} is not restricted to algebra, but also treats the division of an angle and the associated question of the calculation of chords and sines.

In 1579 B\"{u}rgi had entered into the service of the landgrave Wilhelm IV (1532-1592) in Kassel as a watchmaker and instrument maker. After Wilhelm's death he continued his work in the service of Wilhelm's successor Moritz (1572-1632). In 1605 B\"{u}rgi went to Prague and lived there as a watchmaker at the Emperor's Chamber until 1631. Then he returned to Kassel, where he died one year later%
\footnote {The most detailed monograph on B\"{u}rgi's life and work is [Staudacher, 2014].}.

In Kassel, B\"{u}rgi not only worked as an instrument maker but was also involved in astronomical observations and interpretations. As to his mathematical and astronomical knowledge, he seems to have been self-taught, at least he had neither attended a Latin school nor a university. At the time of Wilhelm IV Kassel had become a center of astronomical and mathematical research. Besides B\"{u}rgi also the qualified astronomer and mathematician Christoph Rothmann (\dag 1601)%
\footnote{The date of his death has recently been proved by [Lenke, Roudet, 2014, 233].}
worked at the court. Among the most important temporary visitors we find Paul Wittich (1555? - 1587) and Nicolaus Reimers Ursus (1551-1600), who had both worked with Tycho Brahe in Hven. Wittich, who was in Kassel from 1584 to 1586, brought with him knowledge of \emph{prosthaphaeresis}, a method by which multiplications and divisions can be replaced by additions and subtractions of trigonometrical values. Since it is based on the identity
\begin{equation}
\sin \alpha \cdot \sin \beta = \textstyle\frac{1}{2} \, [\sin \, (90^{\circ} - \alpha + \beta) - \sin \, (90^{\circ} - \alpha - \beta)] \;,\end{equation}
this procedure achieves the same as logarithms, which were invented some decades later%
\footnote{The method of prosthaphaeresis had been introduced by Johannes Werner (1468-1522). For details see [Van Brummelen, 2009, 264-265], with reference to additional secondary literature.}.
During the following years Rothmann and B\"{u}rgi worked with this method. In order to use prosthaphaeresis a sine table is needed, and perhaps this was the reason why B\"{u}rgi studied the calculation of sine values.

Ursus, who lived in Kassel in 1586 and 1587, and B\"{u}rgi became good friends during this period. This can be seen from Ursus's wording in the preface of his \emph{Fundamentum astronomicum}, which was published in 1588 in Strassburg. He writes [Ursus, 1588, fol. *4r]:
``I do not have to explain to which level of comprehensibility this extremely deep and nebulous theory has been corrected and improved by the tireless study of my dear teacher, Justus B\"{u}rgi from Switzerland, by assiduous considerations and daily thought. [...] Therefore neither I nor my dear teacher, the inventor and innovator of this hidden science, will ever regret the trouble and the labor which we have spent.''

The \emph{Fundamentum astronomicum} starts with the mathematical foundations of astronomy. In the first chapter Ursus explains how to count with sexagesimal numbers. In Chapter 2 he deals with the computation of a sine table and in Chapter 3 with the calculation of orthogonal and oblique-angled spherical triangles as well as with plane triangles. The last two chapters are on the observation of the positions of the fixed stars and on the movement of the planets%
\footnote{[Ursus, 1588]. German translation in [Launert, 2012].}.

At the beginning of Chapter 2 Ursus writes [Ursus, 1588, fol. B1v]:
``The calculation of the \emph{Canon Sinuum} can be done [...] either in the usual way, by inscribing the sides of a regular polygon into a circle [...], that is, geometrically. Or it can be done by a special way, by dividing a right angle into as many parts as one wants; and this is arithmetically. This has been found by Justus B\"{u}rgi from Switzerland, the skilful technician of His Serene Highness, the Prince of Hesse.'' However, Ursus does not explain B\"{u}rgi's method in detail, but only makes vague remarks which do not enable the reader to reproduce the procedure. He adds a diagram [Fig. 1], but this does not help much either.


\begin{figure}[h]
\centering
\includegraphics[scale=0.45]{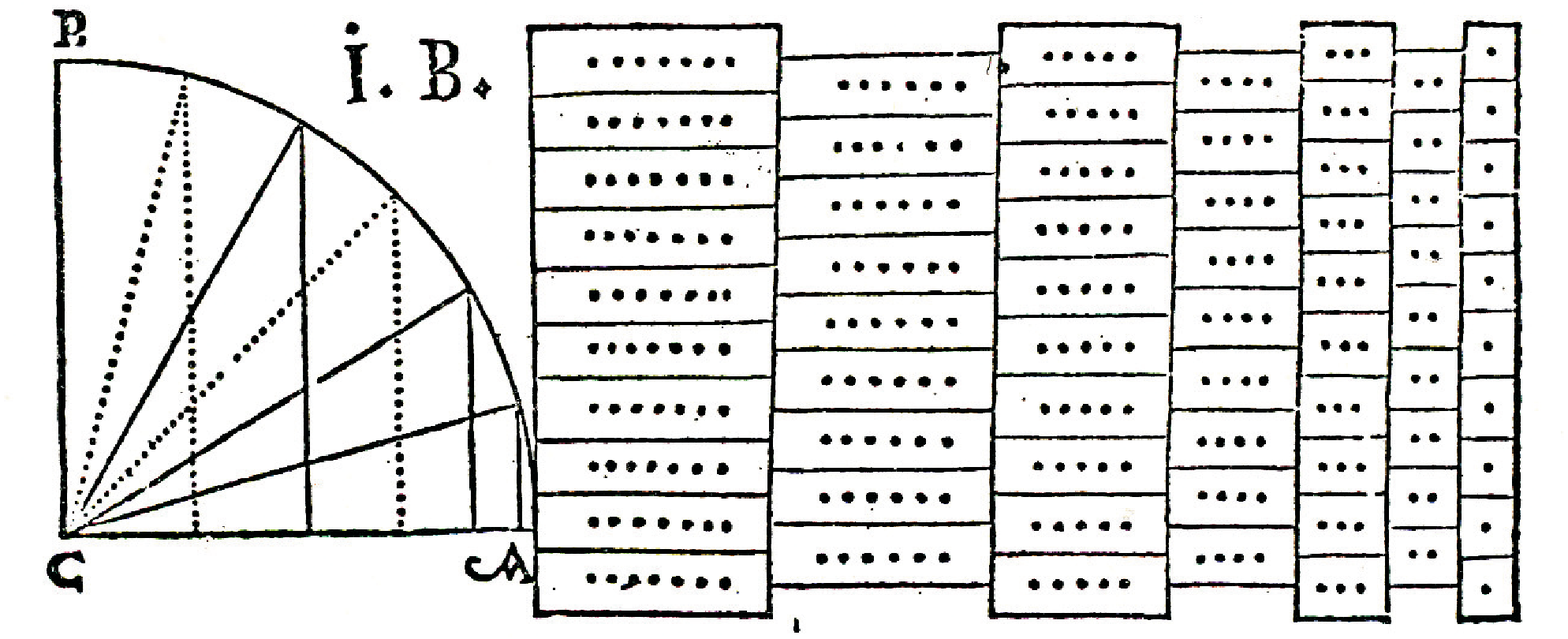}
\caption{Fig. 1. Ursus's diagram (\emph{Fundamentum astronomicum}, fol. C1r)}
\end{figure}

Obviously there was an arrangement between Ursus and B\"{u}rgi that Ursus should only hint to B\"{u}rgi's procedure but not explain it in detail. Based upon a note by B\"{u}rgi himself, nowadays his procedure is usually called his ``Kunstweg'' (``skillful / artful method'')%
\footnote{See [List, Bialas, 1973, 116].}.

B\"{u}rgi never published his method. Ursus's hints motivated Kepler and other scholars to try to reconstruct B\"{u}rgi's procedure, but these attempts were not successful%
\footnote{See [List, Bialas, 1973, 115-122].}.
Modern historians of science have resigned and suggested that it would probably never be possible to decode B\"{u}rgi's ``skillful / artful method''%
\footnote{So [Braunm\"{u}hl, 1900, 210] and [List, Bialas, 1973, 122].}.

Fortunately, this statement is no longer true because a manuscript has survived in which B\"{u}rgi himself explains his method. Menso Folkerts found this hitherto entirely unknown autograph of B\"{u}rgi, with shelfmark IV Qu. 38$^\mathrm{a}$, in the Biblioteka Uniwersytecka in Wroc{\l}aw (Poland). From the preface of the book we know that B\"{u}rgi presented it to the emperor Rudolf II (1552-1612, emperor since 1576) on 22 July 1592 in Prague and that some days later the emperor donated him the amount of 3000 Taler%
\footnote{The details are given in [Folkerts, 2014].}.
Later the manuscript found its way into the library of the Augustinian monastery in Sagan (Lower Silesia; today: \.{Z}aga\'{n}, Poland), the same town where Johannes Kepler lived from 1628 to 1630. In 1810, when the monastery was secularized, the manuscript was brought into the Universit\"{a}tsbibiothek in Breslau (then Germany) where it fell into oblivion.

The manuscript consists of 95 folios and contains a hitherto unknown extensive work on trigonometry by B\"{u}rgi, entitled \emph{Fundamentum Astronomiae}. Already the title reminds of Ursus's book of 1588. Also in the content there is much similarity. Like Ursus, B\"{u}rgi intends to present the basics of astronomical computation, i.e. counting in the sexagesimal system, prosthaphaeresis, the calculation of a sine table and the application of sines to calculate plane and spherical triangles. Unlike Ursus, who wrote in Latin, B\"{u}rgi's work is written in German. A commented edition of B\"{u}rgi's complete text has been published in the ``Abhandlungen'' of  the Bayerische Akademie der Wissenschaften [Launert, 2016].

B\"{u}rgi's treatise consists of two books which are divided into 13 and 11 chapters, respectively. Book 1 deals with logistic numbers, prosthaphaeresis and the calculation of sines. ``Logistic numbers'' are the sexagesimal numbers which were used for astronomical calculations. In the first two chapters B\"{u}rgi explains the four basic operations of arithmetic and the extraction of roots. In the section about multiplication B\"{u}rgi presents on 12 pages a multiplication table in the sexagesimal system from $1\times 1$ to $60\times 60$. The topic of Chapter 3 is prosthaphaeresis. Here B\"{u}rgi proves the central equation (1). The remaining 10 chapters of Book 1 are on sines and especially on the calculation of sine values. Chapters 11 and 12 are the most interesting part of the whole book because B\"{u}rgi explains here in detail his own method for computing in relatively few steps all sine values from 0 to 90 degrees (see below). The result is a sine table for every minute with 5-7 sexagesimal places. Accordingly, this table contains $90 \cdot 60 = 5400$ values and comprises 36 pages, i.e. nearly a fifth of the entire treatise.

Book 2 deals with the calculation of triangles. The first four chapters are on plane triangles and chapters 5-11 on spherical ones. B\"{u}rgi treats the six possible cases of known elements of the triangle.  His purpose was to solve all kinds of spherical triangles and especially triangles whose angles are known. In this he wanted to show that all solutions can be obtained by using only the sine and no other trigonometric function.

Of special interest is the computation of sine values in Book 1. B\"{u}rgi starts in the traditional way: He inscribes regular polygons ($n$ = 3, 4, 5, 6, 10) into a circle and explains how the sines of $\alpha$ = 18$^{\circ}$, 30$^{\circ}$, 36$^{\circ}$, 45$^{\circ}$ and 60$^{\circ}$ can be calculated from these. Then he shows how by means of trigonometric relations the value of $\sin 1\frac{1}{2}^{\circ}$ can be determined. From this B\"urgi finds $\sin 1^{\circ}$, presumably by the traditional method.


\begin{figure}[h]
\centering
\includegraphics[scale=0.45]{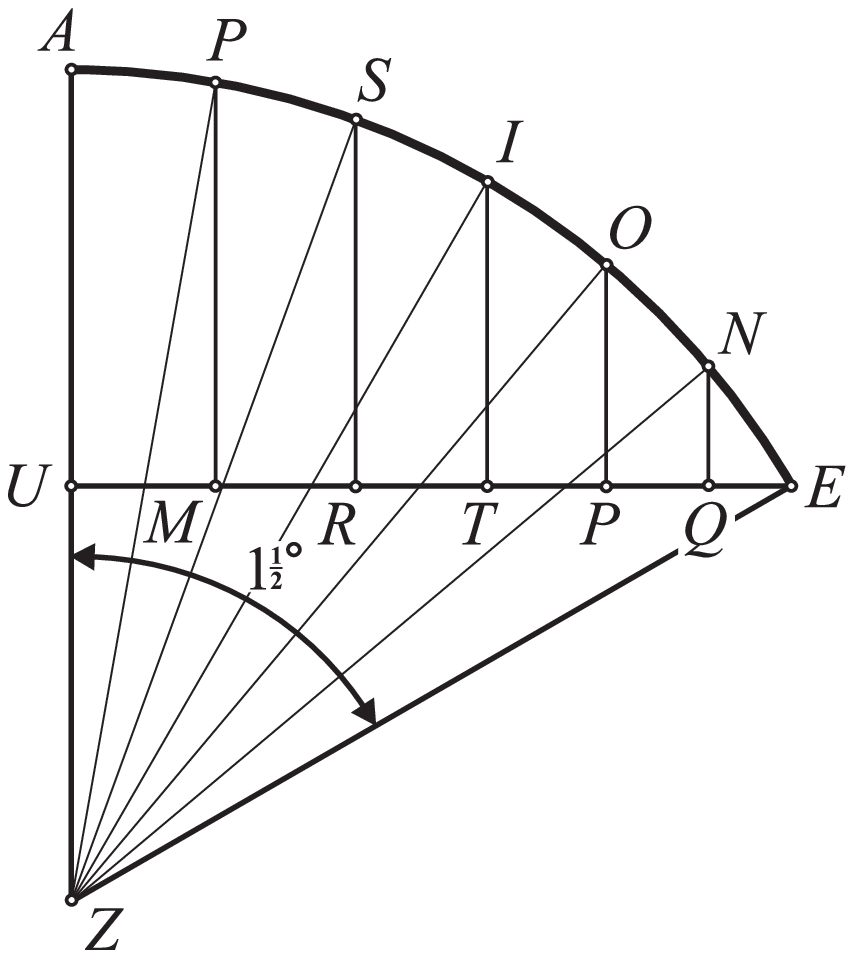}
\caption{Fig. 2. B\"{u}rgi's diagram for calculating $\sin 1^{\circ}$ (\emph{Fundamentum Astronomiae}, fol. 32r)}
\end{figure}

In addition,  in Chapter 10 he presents a new method for finding this value (see Figure 2): He divides the angle of $1 \frac{1}{2}$ degrees, arc \emph{AE} in the figure, into six equal parts of $15'$. If $Z$ is the center of the circle, $U$ the projection of $E$ onto radius $ZA$, and $T$ the projection of $I$, the middle of arc $AE$, onto $UE$, then the line segment $UT$ is equal to $\sin 45'$. B\"{u}rgi finds that $\frac{1}{6}\, UE \approx \sin 15'$ and that $UT + \frac{1}{6}\, UE$ is a good approximation to $\sin 1^{\circ}$.
Using the value of $\sin 1^{\circ}$, B\"{u}rgi is able to calculate in the traditional way -- B\"{u}rgi calls it ``geometrically '' -- a complete sine table.

At the end of Chapter 10, on fol. 34r-v, B\"{u}rgi writes:
``In such a way, with much trouble and labor, the whole \emph{Canon} has been established. For many hundreds of years, up to now, our ancestors have been using this method because they were not able to invent a better one. However, this method is uncertain and dilapidated as well as cumbersome and laborious. Therefore we want to perform this in a different, better, more correct, easier and more cheerful way. And we want to point out now how all sines can be found without the troublesome inscription [of polygons], namely by dividing a right angle into as many parts as one desires.''

B\"{u}rgi describes his ``skilled method'' in Chapter 11 of Book 1. B\"{u}rgi proudly writes that this is ``a very much desired task, which reveals to us the eager thoughts of B\"{u}rgi which have never been seen and which did not exist before''%
\footnote{``postulatum desideratissimum. Sedula quod nobis aperit meditatio Burgi // Quale nec est visum quale nec antefuit'' (fol. 34v).}.

B\"{u}rgi explains his method by an example, namely by the division of a right angle into nine parts. In this way he determines the values of sin 10$^{\circ}$ to sin 90$^{\circ}$. He provides a diagram to show the procedure (see Fig. 3).

\begin{figure}[h]
\centering
\includegraphics[scale=0.5]{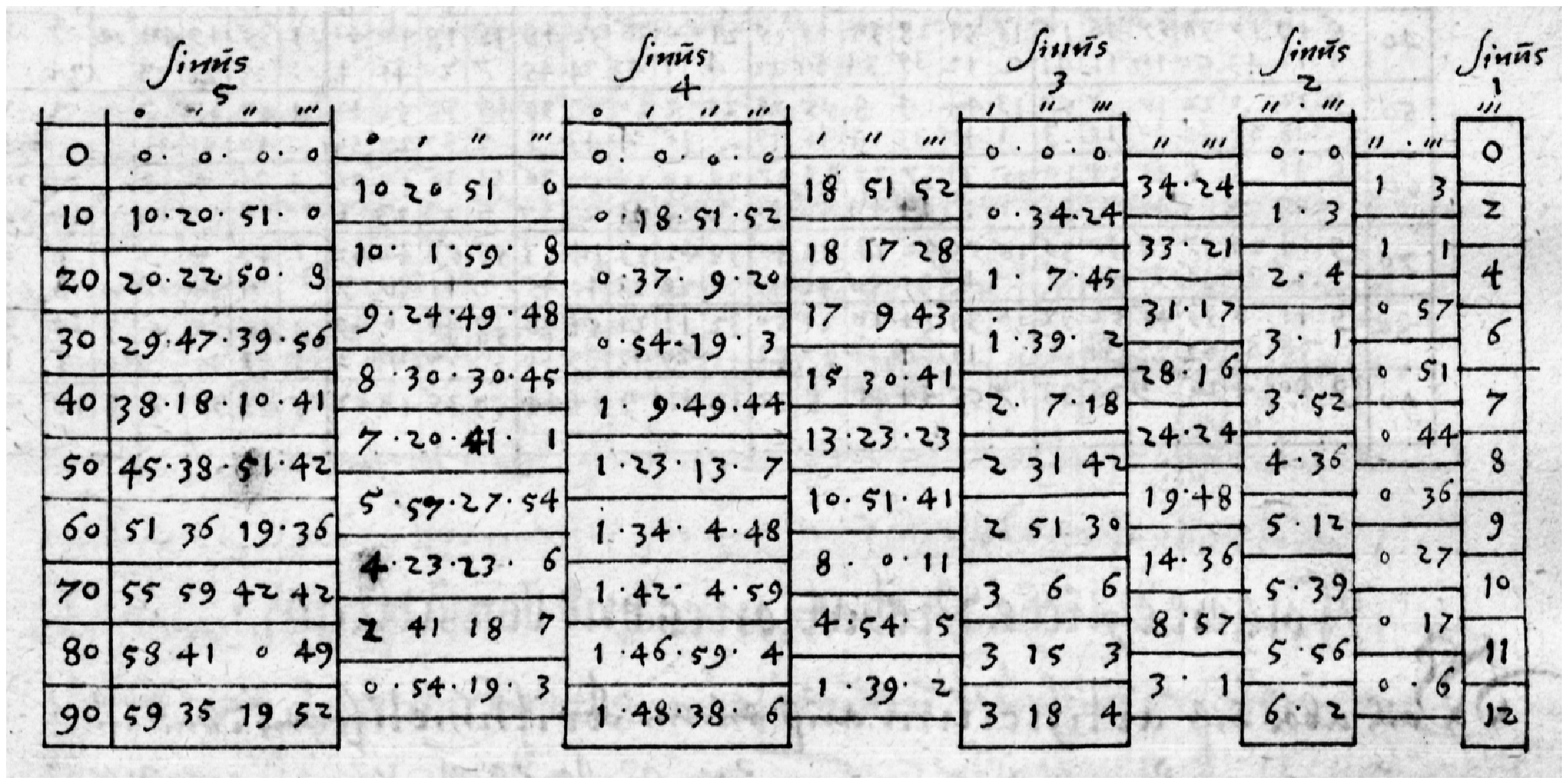}
\caption{B\"{u}rgi's table (\emph{Fundamentum Astronomiae}, fol. 36r)}
\end{figure}


\noindent
B\"{u}rgi writes the numbers in the sexagesimal system. Since we are used to count in the decimal system, the following notes refer to [Fig. 4].

{\footnotesize$$\begin{array}{ccccccccc}
{\addtolength{\arraycolsep}{-0.550em}\begin{array}{r|r|}
\cline{2-2}&\multicolumn{1}{c|}{\text{Column}\;5}\\
\hline\vline\hspace{0.8em}0\;&0\:\\
\hline\vline\;10\;&2,\!235,\!060\:\\
\hline\vline\;20\;&4,\!402,\!208\:\\
\hline\vline\;30\;&6,\!435,\!596\:\\
\hline\vline\;40\;&8,\!273,\!441\:\\
\hline\vline\;50\;&9,\!859,\!902\:\\
\hline\vline\;60\;&\:11,\!146,\!776\:\\
\hline\vline\;70\;&12,\!094,\!962\:\\
\hline\vline\;80\;&12,\!675,\!649\:\\
\hline\vline\;90\;&12,\!871,\!192\:\\
\hline
\end{array}}
{\addtolength{\arraycolsep}{-0.550em}\begin{array}{r}
\\
\hline\:\:2,\!235,\!060\:\\
\hline\:2,\!167,\!148\:\\
\hline\:2,\!033,\!388\:\\
\hline\:1,\!837,\!845\:\\
\hline\:1,\!586,\!461\:\\
\hline\:1,\!286,\!874\:\\
\hline\:948,\!186\:\\
\hline\:580,\!687\:\\
\hline\:195,\!543\:\\
\hline\end{array}}
{\addtolength{\arraycolsep}{-0.550em}\begin{array}{|r|}
\hline\multicolumn{1}{|c|}{\text{Col.}\;4}\\
\hline\:0\:\\
\hline\:67,\!912\:\\
\hline\:\:133,\!760\:\\
\hline\:195,\!543\:\\
\hline\:251,\!384\:\\
\hline\:299,\!587\:\\
\hline\:338,\!688\:\\
\hline\:367,\!499\:\\
\hline\:385,\!144\:\\
\hline\:391,\!086\:\\
\hline
\end{array}}
{\addtolength{\arraycolsep}{-0.550em}\begin{array}{r}
\\
\hline\:\:67,\!912\:\\
\hline\:65,\!848\:\\
\hline\:61,\!783\:\\
\hline\:55,\!841\:\\
\hline\:48,\!203\:\\
\hline\:39,\!101\:\\
\hline\:28,\!811\:\\
\hline\:17,\!645\:\\
\hline\:5,\!942\:\\
\hline\end{array}}
{\addtolength{\arraycolsep}{-0.550em}\begin{array}{|r|}
\hline\multicolumn{1}{|c|}{\text{Col.}\;3}\\
\hline\:0\:\\
\hline\:2,\!064\:\\
\hline\:4,\!065\:\\
\hline\:5,\!942\:\\
\hline\:7,\!638\:\\
\hline\:9,\!102\:\\
\hline\:\:10,\!290\:\\
\hline\:11,\!166\:\\
\hline\:11,\!703\:\\
\hline\:11,\!884\:\\
\hline
\end{array}}
{\addtolength{\arraycolsep}{-0.550em}\begin{array}{r}
\\
\hline\:\:2,\!064\:\\
\hline\:2,\!001\:\\
\hline\:1,\!877\:\\
\hline\:1,\!696\:\\
\hline\:1,\!464\:\\
\hline\:1,\!188\:\\
\hline\:876\:\\
\hline\:537\:\\
\hline\:181\:\\
\hline\end{array}}
{\addtolength{\arraycolsep}{-0.550em}\begin{array}{|r|}
\hline\multicolumn{1}{|c|}{\:\text{Col.}\,\!2\:}\\
\hline\:0\:\\
\hline\:63\:\\
\hline\:\:124\:\\
\hline\:181\:\\
\hline\:232\:\\
\hline\:276\:\\
\hline\:312\:\\
\hline\:339\:\\
\hline\:356\:\\
\hline\:362\:\\
\hline
\end{array}}
{\addtolength{\arraycolsep}{-0.550em}\begin{array}{r}
\\
\hline\:\:63\:\\
\hline\:61\:\\
\hline\:57\:\\
\hline\:51\:\\
\hline\:44\:\\
\hline\:36\:\\
\hline\:27\:\\
\hline\:17\:\\
\hline\:6\:\\
\hline\end{array}}
{\addtolength{\arraycolsep}{-0.550em}\begin{array}{|r|}
\hline\:\text{Col.}\,\!1\:\\
\hline\:0\:\\
\hline\:2\:\\
\hline\:4\:\\
\hline\:6\:\\
\hline\:7\:\\
\hline\:8\:\\
\hline\:9\:\\
\hline\:10\:\\
\hline\:11\:\\
\hline\:12\:\\
\hline
\end{array}}
\end{array}$$}

\begin{center}\footnotesize Figure 4: B\"{u}rgi's table (decimal system)\end{center}

\noindent
From a modern point of view B\"{u}rgi describes an iteration method which starts with 10 arbitrary numbers $a_{1i} $ (with $a_{10} = 0$ and $a_{11}$, $a_{12}, \ldots , a_{19} \in \mathbb {N}$) and leads to approximations for sin 10$^{\circ}$, \ldots , sin 90$^{\circ}$. B\"{u}rgi's method is a strongly convergent process for finding the sine values, but B\"{u}rgi does not prove why it works and his formulation is not easy to understand%
\footnote{Dieter Launert has decoded B\"{u}rgi's text as explained below and Andreas Thom has found a modern proof which is given at the end of this paper.}.

B\"{u}rgi presents a series of approximations, whose intermediate results are written in the columns labeled 1 to 5 in Fig. 3. Between each two consecutive approximations an auxiliary column is inserted. Using only additions and halving, B\"{u}rgi proceeds in the following way:

1) In the first column on the right, labeled ``Column 1'', a zero is placed in the upper cell. The other nine cells are filled with an arbitrary series of natural numbers, for which B\"{u}rgi selects 2, 4, 6, 7, 8, 9, 10, 11, 12. B\"{u}rgi does not mention that the numbers in column 1 are in principle his first approximations for sin 0$^{\circ}$ to sin 90$^{\circ}$ in relation to the value 12 for sin 90$^{\circ}$. Nor does he mention that the calculation will be shorter if the numbers correspond approximately to the ratio of the sine values of 10$^{\circ}$ to 90$^{\circ}$.

2) The number in the bottom cell of column 1 is halved, and this halve is placed in the bottom of a new, auxiliary column to the left of column 1. This column is displaced upwards by half a row, so that in B\"{u}rgi's example the 6 (halve of 12) in the lowest cell is placed in between the 11 and 12 in the last two cells of column 1.

3) The auxiliary column thus begun, which has no heading, is now filled step by step from bottom to top by adding in each case the previously found number in the auxiliary column to the number in column 1 placed just above it. In B\"{u}rgi's example we receive the numbers: $6+11 = 17$, $17+10 = 27$, $27+9 = 36$, $36+8 = 44$, $44+7 = 51$, $51+6 = 57$, $57+4 = 61$, $61+2 = 63$. These values correspond, respectively, to estimates for the sines of the intermediate angles $5^{\circ}$, $15^{\circ}$, \ldots $85^{\circ}$. Again, B\"{u}rgi does not mention this explicitly.

4) Now a new results column, headed ``Column 2'', is added to the left of the first auxiliary column. A zero is placed in the upper cell of this column, displaced slightly upwards from the upper number in the auxiliary column. Then, now from top to bottom, the cells in column 2 are filled step by step with the sum of the previously found number and the number in the auxiliary column just below it. Thus we have: $0+63 = 63$, $63+61 = 124$, $124+57 = 181$, etc., up to $356+6 = 362$. The numbers in column 2 are a second approximation, better than the initial one in column 1, to the sine values for 10$^{\circ}$, 20$^{\circ}$, \ldots , 80$^{\circ}$, 90$^{\circ}$ in relation to the value 362 for the sine of 90$^{\circ}$.

5) This procedure is now repeated: In order to obtain a third approximation, half of the lowest number in column 2, namely number 181, is entered into the lowest cell of a new auxiliary column to the left of column 2, which is again displaced by half a row. Then, as in step 3), first the ascending sums are formed and entered into the new auxiliary column, leading to 2,064 in the top cell. Thereafter column 3 is built with a zero in the top cell, and it is filled by adding to the respectively found numbers the numbers just below them in the auxiliary column. Thus this column starts with $0 + 2,\!064 = 2,\!064$ and ends with $11,\!703 + 181 = 11,\!884$. Now column 3 contains a third approximation to the 9 sine values in proportion to the highest number 11,884.

6) B\"{u}rgi emphasizes that the numbers in columns 1, 2, 3, \dots do not directly give the sine values. Instead, they are in relation to the ``radius'', i.e. the number given in the last line, which corresponds to the sine of 90$^{\circ}$. Thus, the values in each column must be divided by the number in the bottom cell.

7) This procedure can be repeated any number of times in order to determine the sines to a higher accuracy. In B\"{u}rgi's example one obtains in column 5, after the 4th iteration:

sin 10$^{\circ}$ = 2,235,060 : 12,871,192 = 0.17364825 (exact: 0.17364818\dots)

.....

sin 80$^{\circ}$ = 12,675,649 : 12,871,192 = 0.984807701 (exact: 0.984807753\dots).

\noindent
Thus the obtained values are correct to 6 or 7 places. An accuracy of 9-10 decimal places can be reached, if one continues the procedure up to column 8. This is still a reasonable computational effort.

In this way B\"{u}rgi has found an arithmetic procedure for computing sine values with arbitrary accuracy. By dividing the right angle into 90 parts, B\"{u}rgi is able to calculate the sines of all degrees from sin 1$^{\circ}$ to sin 90$^{\circ}$.

In Chapter 11 B\"{u}rgi also deals with the question of how to calculate sine values for every minute%
\footnote{On f.38r-40v.}.
If he tried to do this by his method, he would have to divide the right angle into 5400 parts, making the calculations very tedious. Instead, B\"{u}rgi finds a type of procedure in order to compute the sine of 1 minute as exactly as possible. He divides the known value for sin 1$^{\circ}$ by 60 to receive an approximation for sin 1$'$. This he improves by means of two corrections and obtains a sufficiently exact value (sin 1$'$ = 0.000,290,888,63 the exact 0.000,290,888,20). With the help of trigonometric relations he then computes sin 2$'$, etc. However, B\"{u}rgi shows that it is not necessary to continue this procedure up to sin 59$'$: by using first and higher order differences of consecutive sine values he derives a simple relation for producing the further sines of minutes of arc\footnote{For details, see the commented edition of B\"{u}rgi's text in [Launert, 2016, 63-66].}. The result is B\"{u}rgi's sine table for every minute, i.e. with $90 \cdot 60 = 5400$ entries, which is extant in his \emph{Fundamentum Astronomiae}.

Because the sine table in Rheticus's \emph{Opus Palatinum}, which had been printed in 1596, contained many errors, B\"{u}rgi later computed another, even more detailed sine table for every two seconds of arc, thus containing 162,000 sine values%
\footnote{See [List, Bialas, 1973, 112f.]. B\"{u}rgi mentions his reason for calculating this table in his \emph{Coss} (see [List, Bialas, 1973, 7]).}.
This table has not survived.

B\"{u}rgi took special care that his method of calculating sines did not become public in his time. Only since his \emph{Fundamentum Astronomiae} has been found, do we know with certainty that he made use of tabular differences and are we acquainted with the details of his procedure. But B\"{u}rgi's method was not entirely unknown: Dieter Launert has found evidence that a trace of his work leads to England and to the mathematician Henry Briggs (1561-1630), who is well known for his role in the history of logarithms: A copy of Raimarus Ursus's \emph{Fundamentum astronomicum}, which is kept in the University Library of Leiden (The Netherlands), includes a folio, written about 1620, with numbers and calculations, a table similar to B\"{u}rgi's diagram for calculating the sines (Fig. 3), and a name which is most likely ``H. Briggs''. This proves that Briggs was acquainted with B\"{u}rgi's method. Very likely Briggs obtained this information from John Dee (1527-1608), who was in contact with Rothmann, when he visited Kassel in 1586 and 1589\footnote{For details, see [Launert, 2016, ...].}. It seems possible that the method of differences that Briggs used for the calculation of logarithms%
\footnote{See [Sonar, 2002, 197-201].} has some relation to B\"{u}rgi's method of computing sines.

In his \emph{Fundamentum Astronomiae} B\"{u}rgi does not give any indication of how he had found his ``skillful / artful method'' and why it works. A modern proof of its correctness, provided by Andreas Thom, follows below%
\footnote{Peter Ullrich (Koblenz) has made a conjecture how B\"{u}rgi could have found his method, which he intends to publish elsewhere.}.


\section{Modern proof}

Let $n$ be a natural number and $(a_1,\dots,a_n)$ be a sequence (or vector) of non-negative real numbers. Let's assume that there exists some $1 \leq j \leq n$ such that $a_j \neq 0$. B\"urgi describes an algorithm that associates with the sequence $(a_1,\dots,a_n)$ a new sequence $(a'_1,\dots,a'_n)$. His computation of the sine-function at points $j \pi/(2n)$ for $1 \leq j \leq n$ relies on the fact that the iterated application of this assignment yields a sequence of vectors that -- after a suitable normalization of the last entry of the vector -- converges quickly to the vector $$\left(\sin\left( \frac{\pi}{2n}\right),\sin\left( \frac{\pi}{n}\right),\sin\left( \frac{3\pi}{2n}\right),\dots,\sin\left( \frac{(n-1)\pi}{2n}\right),1\right).$$
Let's first describe the algorithm in mathematical terms. The algorithm proceeds in two steps.

\begin{enumerate}
\item Divide $a_n$ by two and set $b_n:=a_n/2$. Then, in an descending recursion define $b_j := a_j + b_{j+1}$ for $1 \leq j \leq n-1$. This way, one obtains a sequence $(b_1,\dots,b_n)$ with
$$b_j = \sum_{k=j}^{n-1} a_k + \frac{a_n}2.$$
\item Set $a'_1 :=b_1$ and define in an ascending recursion $a'_j := b_j + a'_{j-1}$ for all $2 \leq j \leq n$. We obtain a sequence $(a'_1,\dots,a'_n)$, which is defined as
$$a'_j := \sum_{l=1}^j b_l.$$
\end{enumerate}
Taking the two steps together, we obtain
$$a'_j = \sum_{l=1}^j b_l = \sum_{l=1}^j \sum_{k=j}^{n-1} a_k + \frac{j \cdot a_n}2 = \sum_{l=1}^{j-1} l \cdot a_l + \sum_{l=j}^{n-1} j \cdot a_l + \frac{j \cdot a_n}2.$$
In other words, there is a linear map  $(a_1,\dots,a_n) \mapsto (a'_1,\dots,a'_n)$, which is given by left-multiplication with an $n\times n$-matrix of the form
$$M := \left(\begin{matrix}
1 & 1 & 1 & 1 &\cdots & 1 & 1/2 \\
1 & 2 & 2 & 2 &\cdots & 2 & 1 \\
1 & 2 & 3 & 3 &\cdots & 3 & 3/2 \\
1 & 2 & 3 & 4 &\cdots & 4 & 2 \\
\vdots & \vdots &\vdots &\vdots & \ddots & \vdots& \vdots& \\
1 & 2 & 3 & 4 & \cdots & n-1 & (n-1)/2 \\
1 & 2 & 3 & 4 & \cdots & n-1 & n/2
\end{matrix} \right).$$
That is, we obtain $(a'_1,\dots,a'_n)^t = M(a_1,\dots,a_n)^t$. Since the matrix $M$ has only positive entries, the Theorem of Perron-Frobenius implies that the iterated application of the mapping $M$ leads to a sequence of vectors
$$(a^{(k)}_1,\dots,a^{(k)}_n)^t := M^k(a_1,\dots,a_n)^t,$$
such that the ratios $a^{(k)}_j/a^{(k)}_n$ converge to positive real numbers $\lambda_j$ as $k$ tends to infinity. Moreover, the vector $(\lambda_1,\dots,\lambda_{n-1},1)$ is an eigenvector of $M$ with respect to the unique eigenvalue of maximal modulus -- the so-called Perron-Frobenius eigenvalue.

B\"urgi's method relies on the fact that $\lambda_j = \sin\left( j \pi/(2n)\right)$ for all $1 \leq j \leq n$.
In order to see this, we will show that the vector
$$v:=\left(\sin\left( \frac{\pi}{2n}\right),\sin\left( \frac{\pi}{n}\right),\sin\left( \frac{3\pi}{2n}\right),\dots,\sin\left( \frac{(n-1)\pi}{2n}\right),1\right)$$ is an eigenvector of $M$ ist. If that is shown, it will imply that $v$ is indeed an eigenvector with respect to Perron-Frobenius eigenvalue, since eigenvectors with respect to any other eigenvalue must have a non-positive or non-real entry. Thus, since the eigenspace with respect to the Perron-Frobenius eigenvalue is one-dimensional, we obtain
$$\lambda_j = \sin\left( \frac{j\pi}{2n}\right)$$ for all $1 \leq j \leq n$ and the proof is complete.

Looking at the last coefficient of the vector $Mv^t$, we expect that the Perron-Frobenius eigenvalue is equal to
\begin{equation}\label{eq:eigenwert}
  \mu := \sum_{l=1}^{n-1} l \cdot \sin\left( \frac{l \pi}{2n}\right) + \frac n2 = \frac{\csc^2(\pi/(4n))}4>0.
\end{equation}

In order to prove the second part of Equation \eqref{eq:eigenwert} and also in order to complete the other computations, we will need the following observations concerning telescoping sums. For real numbers $a_0,\ldots,a_{n}$, where $a_0=0$ and $1 \leq j \leq n$, we have:\begin{equation}\label{eq:teleskop1}
   \sum_{l=1}^{j-1} l \cdot a_{l}- \sum_{l=1}^{j-1}\frac{l}{2} \cdot (a_{l-1}+a_{l+1})=\frac{1}{2}(j \cdot a_{j-1}-(j-1) \cdot a_j),
\end{equation}
\begin{equation}\label{eq:teleskop2}
\sum_{l=j}^{n-1}a_l- \sum_{l=j}^{n-1}\frac12 \left(a_{l-1}+ a_{l+1} \right) =\frac 12\left(a_{n-1}-a_n+a_j-a_{j-1}\right),\end{equation}
results that can be proved without any problems by induction.

It follows from the well-known addition theorem $$\sin(x)+\sin(y)=2\sin\left(\frac{x+y}2\right)\cos\left(\frac{x-y}2\right)$$ for real $x,y$ that:
\begin{equation}\label{eq:cos_mittel}
\cos\left( \frac{\pi}{2n}\right) \cdot \sin\left( \frac{l \pi}{2n}\right) = \frac12 \left(\sin\left( \frac{(l-1) \pi}{2n}\right) + \sin\left( \frac{(l+1) \pi}{2n}\right) \right)
\end{equation}
for all $1 \leq l \leq n$.
Since $2 \cdot \sin^2(\pi/(4n)) = 1 - \cos\left( \frac{\pi}{2n}\right)$, the second part of Equation \eqref{eq:eigenwert} is implied by the following computation.
\begin{eqnarray*}
&&2\sin^2\left(\frac{\pi}{4n}\right) \cdot  \mu \\& =& \left(1 - \cos\left( \frac{\pi}{2n}\right)\right) \cdot \mu \\
  & =& \sum_{l=1}^{n-1} l \cdot \sin\left( \frac{l \pi}{2n}\right) + \frac n2- \cos\left( \frac{\pi}{2n}\right)\left(\sum_{l=1}^{n-1} l \cdot \sin\left( \frac{l \pi}{2n}\right) + \frac n2\right)\\
  & \stackrel{\eqref{eq:cos_mittel}}{=}&\sum_{l=1}^{n-1} l \cdot \sin\left( \frac{l \pi}{2n}\right)
   \\ &&\quad - \sum_{l=1}^{n-1} l \cdot \frac12 \left(\sin\left( \frac{(l-1) \pi}{2n}\right) + \sin\left( \frac{(l+1) \pi}{2n}\right) \right) \\
   && \quad +\frac n2 - \frac n2 \cdot \cos\left( \frac{\pi}{2n}\right)\\
  &\stackrel{\eqref{eq:teleskop1}}{=}& \frac{1}{2} \left(n \cdot \sin\left( \frac{(n-1) \pi}{2n} \right)- (n-1) \cdot \sin\left( \frac{n \pi}{2n} \right)\right) \\
  && \quad +\frac{n}{2} - \frac n2 \cdot \cos\left( \frac{\pi}{2n}\right)
  \\&=&\frac12.
\end{eqnarray*}

Now, in order to show that $v$ is indeed an eigenvector of $M$ with eigenvalue $\mu$, we need to show that for all $1 \leq j \leq n-1$, the following equation holds.
\begin{eqnarray} \label{eqeig}
\mu \cdot \sin\left( \frac{j \pi}{2n}\right)
&=& \sum_{l=1}^{j-1} l \cdot \sin\left( \frac{l \pi}{2n}\right) + \sum_{l=j}^{n-1} j \cdot \sin\left( \frac{l \pi}{2n}\right) + \frac{j}2.
\end{eqnarray}
After multiplication with $2 \sin^2(\pi/(4n)) = 1 - \cos\left( \frac{\pi}{2n}\right)$ we obtain the equivalent equation:
\begin{align*}
\frac12 \sin\left( \frac{j \pi}{2n}\right) &=
\left(1 - \cos\left( \frac{\pi}{2n}\right) \right) \cdot \left(\sum_{l=1}^{j-1} l \cdot \sin\left( \frac{l \pi}{2n}\right) + \sum_{l=j}^{n-1} j \cdot \sin\left( \frac{l \pi}{2n}\right) + \frac{j}2 \right).
\end{align*}
In order to prove this equivalent statement, we compute as follows.
\begin{eqnarray*}
&&\left(1 - \cos\left( \frac{\pi}{2n}\right) \right) \cdot \left(\sum_{l=1}^{j-1} l \cdot \sin\left( \frac{l \pi}{2n}\right) + \sum_{l=j}^{n-1} j \cdot \sin\left( \frac{l \pi}{2n}\right) + \frac{j}2 \right)\\
\\ &\stackrel{\eqref{eq:cos_mittel}}{=}& \sum_{l=1}^{j-1} l \cdot \sin\left( \frac{l \pi}{2n}\right) + \sum_{l=j}^{n-1} j \cdot \sin\left( \frac{l \pi}{2n}\right) + \frac{j}2
  \\ &&\quad -\sum_{l=1}^{j-1} l \cdot \frac12 \left(\sin\left( \frac{(l-1) \pi}{2n}\right) + \sin\left( \frac{(l+1) \pi}{2n}\right) \right)
  \\ && \quad- \sum_{l=j}^{n-1} j \cdot \frac12 \left(\sin\left( \frac{(l-1) \pi}{2n}\right)+ \sin\left( \frac{(l+1) \pi}{2n}\right) \right) - \cos\left( \frac{\pi}{2n}\right)\frac{j}2
\end{eqnarray*}
Applying now Equations \eqref{eq:teleskop1} and \eqref{eq:teleskop2}, we obtain:
\begin{eqnarray*}
&=& \frac12\left(j\cdot \sin\left( \frac{(j-1) \pi}{2n}\right)-(j-1)\cdot \sin\left( \frac{j \pi}{2n}\right) \right)
 \\ && \quad  + \frac{j}2\left(\sin\left( \frac{(n-1) \pi}{2n}\right)-\sin\left( \frac{n \pi}{2n}\right)+\sin\left( \frac{j \pi}{2n}\right)-\sin\left( \frac{(j-1) \pi}{2n}\right)\right)
 \\ && \quad + \frac{j}2 -\cos\left( \frac{\pi}{2n}\right)\cdot \frac{j}2
  \\ & =& \frac12 \sin\left( \frac{j\pi}{2n}\right).
\end{eqnarray*}
This proves Equation \eqref{eqeig} and hence, that $v^t$ is an eigenvector of $M$. We can now conclude that B\"urgi's algorithm provides indeed an iterative procedure to compute the sine-function at points $j\pi/(2n)$ for $1 \leq j \leq n$ to arbitrary precision.

It is natural to wonder how B\"urgi might have found this algorithm -- and we think that our proof of its correctness sheds some light on this problem too. We thank one unknown referee of this paper for pointing out to us that the inverse of the matrix $M$ has the following form
$$M^{-1} = \left( \begin{matrix} 
2 & -1 & 0 & 0 & \hdots & 0 & 0 & 0\\
-1 & 2 & -1 & 0&\hdots& 0  & 0 & 0\\
0 & -1 & 2 &  -1& \ddots& 0 & 0 & 0 \\
0 & 0 & -1 &2 & \ddots& \ddots & 0  & 0\\
 \vdots & \vdots& \ddots & \ddots & \ddots& \ddots & \ddots & \vdots\\
 0 & 0 & 0 & \ddots & \ddots &2 & -1 & 0\\
 0 & 0 &0 & 0 &\ddots & -1 &  2& -1\\
 0 & 0 &0 & 0 &\hdots & 0 & -2 &2\\
\end{matrix}\right).$$
Thus, the matrix $M^{-1}$ computes up to a sign and some normalization at the start and end of the sequence just second differences, a process that was well-known to B\"urgi and also well-known in the study of sine-tables. As a consequence, we see that B\"urgi's algorithm just reverses the process of forming second differences, i.e., performs up to sign some form of two-fold discrete integration -- with the right normalization at the start and end of the sequence. Of course, our Perron-Frobenius eigenvector $v$ of $M$ is also an eigenvector of $M^{-1}$, but now for the {\it smallest} eigenvalue. B\"urgi's insight must have been that the study of iterations of $M$ is much more useful than those of $M^{-1}$ in order to approximate the entries of the critical eigenvector. Indeed, as we have seen, this process has unexpected stability properties leading to a quickly convergent sequence of vectors that approximate this eigenvector and hence the sine-values. The reason for its convergence is more subtle than just some geometric principle such as exhaustion, monotonicity, or Newton's method, it rather relies on the equidistribution of a diffusion process over time  -- an idea which was later formalized as the Perron-Frobenius theorem and studied in great detail in the theory of Markov chains. Thus, we must admit that B\"urgi's insight anticipates some aspects of ideas and developments that came to full light only at the beginning of the 20th century.
\section*{Acknowledgments}

We would like to thank Benno van Dalen (Munich), Stefan Deschauer (Dresden), Peter Ullrich (Koblenz) and Marcus Waurick (Dresden) for their help.

\section*{Bibliography}

Braunm\"{u}hl, A. v., 1900. Vorlesungen \"{u}ber Geschichte der Trigonometrie. Erster Teil. Von den \"{a}ltesten Zeiten bis zur Erfindung der Logarithmen. Teubner, Leipzig.

Curtze, M., 1900. Urkunden zur Geschichte der Trigonometrie im christ\-li\-chen Mittelalter. In: Bibliotheca Mathematica, 3. Folge, 1, 321-416.

Folkerts, M., 2014. Eine bisher unbekannte Schrift von Jost B\"{u}rgi zur Trigonometrie. In: Gebhardt, R. (Ed.), Arithmetik, Geometrie und Algebra in der fr\"{u}hen Neuzeit. Adam-Ries-Bund, Annaberg-Buchholz, pp.107-114.

Launert, D., 2012. Astronomischer Grund. Fundamentum Astro\-no\-mi\-cum (1588) des Nicolaus Reimers Ursus. Deutsch, Frankfurt am Main.

Launert, D., 2016. B\"{u}rgis Kunstweg im Fundamentum Astronomiae. Entschl\"{u}sselung eines R\"{a}tsels. Verlag der Bayerische Akademie der Wissenschaften, M\"{u}nchen.

Lenke, N., Roudet, N., 2014. Johannes, Christoph und Bartholomaeus Rothmann. Einige biographische Erg\"{a}nzungen zu den Gebr\"{u}dern Rothmann. In: Dick, W.R., F\"{u}rst, D. (Eds.), Lebensl\"{a}ufe und Himmelsbahnen. Fest\-schrift zum 60. Geburtstag von J\"{u}rgen Hamel. Akademische Ver\-lags\-an\-stalt, Leipzig, pp. 223-242.

List, M., Bialas, V., 1973. Die Coss von Jost B\"{u}rgi in der Redaktion von Johannes Kepler. Verlag der Bayerischen Akademie der Wissenschaften, M\"{u}nchen.

North, J.D., 1976. Richard of Wallingford. An edition of his writings with introductions, English translation and commentary. 3 vols. Clarendon Press, Oxford.

Pedersen, F.S., 2002. The Toledan Tables. A review of the manuscripts and the textual versions with an edition. 4 volumes. Reitzels, Copenhagen.

Sonar, Th., 2002. Der fromme Tafelmacher. Die fr\"{u}hen Arbeiten des Henry Briggs. Logos Verlag, Berlin.

Staudacher, S., 2014. Jost B\"{u}rgi, Kepler und der Kaiser. Verlag Neue Z\"{u}rcher Zeitung, Z\"{u}rich.

Toomer, G.J., 1968. A survey of the Toledan Tables. In: Osiris 15, 5-174.

Ursus, N.R., 1588. Fundamentum Astronomicum. Jobin, Stra{\ss}burg 1588.

Van Brummelen, G., 2009. The Mathematics of the Heavens and the Earth. The Early History of Trigonometry. Princeton University Press, Princeton, Oxford.

\end{document}